\documentclass[nobib]{tufte-handout}
\usepackage{amsfonts,amssymb,amsmath,amscd}
\usepackage{graphicx}
\usepackage[backend=biber]{biblatex}
\date{}

\newcommand{\acknowledgmark}{\marginnote{The author thanks Claude Opus 4 and 4.6 (Anthropic) for extensive discussions that helped develop and refine these ideas.}}

\begin{document}

\begin{fullwidth}
\noindent{\huge\itshape On the Mechanical Creation of Mathematical Concepts}

\vspace{0.5em}
\noindent{\large Asvin G.}\marginnote{Institute for Advanced Study\\ \texttt{gasvinseeker94@gmail.com}}
\end{fullwidth}

\vspace{3.5em}
\acknowledgmark

The formalization of intelligence has a long history starting at least with Leibniz and passing through Boole, Frege, Hilbert and Turing amongst many others. The 19th century culminated with Hilbert's program attempting to formalize all of mathematics, and thus give it certainty - a dream that was destroyed in large part by G\"odel's incompleteness theorems. The incompleteness theorems however were the seed for the development of computability theory in the twentieth century, which has led to a practical mechanization of large parts of "intelligent work". In the twenty-first century, however, we have seen a shift from mechanizing the execution of such work to mechanizing the learning process itself. Taking mathematics as a case study, we saw the automation of the verification of mathematical proof in the twentieth century, and now we are seeing the automation of the discovery of mathematical proof in the twenty-first. To what extent will the current procedures lead to a complete automation of mathematical discovery? Answering this will first require us to answer what the current human process of mathematical discovery consists of, and where AI systems today stand in relation to this schema. This essay will try to answer both questions, in an attempt to speak to questions of intelligence more fundamentally.

\section*{Priors, search, and information}

Any act of problem-solving, from a game of chess to a mathematical proof, combines two things: the knowledge the solver brings to the problem beforehand, and the computation performed on the specific instance at hand. I will call these \emph{priors} and \emph{local search}.
\begin{marginfigure}[0cm]
\centering
\includegraphics[width=\linewidth]{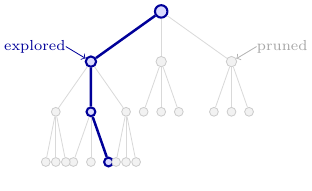}
\caption{Priors and search. Priors determine which branches look promising (dark) and which to skip (light). Local search explores the promising branches, finding an outcome (\checkmark) without traversing the full tree.}
\end{marginfigure} These two are familiar, but there is a third element that will concern me more: the process by which the solver extracts information from search to update understanding.

In chess, a strong player arrives at each position with a repertoire of patterns learned through years of study: typical pawn structures, dangerous piece configurations, favorable opening sequences. This is prior knowledge. But priors alone do not determine the move, since each position is specific and the player must calculate concrete variations: if I take that bishop, what follows after ...Nxe4, Bxf7+, Kd8? This calculation is local search. The two work together: heuristics encode experience into rapid pattern recognition but can miss features of a particular position, while explicit calculation is slower and error-prone but can uncover possibilities that heuristics overlook, and they feed into each other, with heuristics selecting which lines are worth calculating and the results of calculation sharpening the heuristics.

AlphaGo builds this dichotomy explicitly into its architecture \cite{silver2017mastering}. A neural network trained on millions of games learns a \emph{policy}, a probability distribution over moves in any position, which encodes trained intuition about which moves look promising. Monte Carlo tree search then explores continuations from the current position, playing out hypothetical futures. During a single game, the tree search uses the policy to guide local computation but does not modify the policy itself, so the concepts embedded in the network stay fixed while search unfolds. Between games, however, the outcomes of search feed back into the policy through training, and over millions of rounds of self-play this reshapes the system's priors: search produces outcomes, outcomes update the policy, the updated policy guides better search in the next game.

Human chess players work similarly. A strong player's concepts (``bad bishop,'' ``color complex,'' ``overextended pawn chain'') are learned over hundreds or thousands of games. During any single game, these concepts are fixed, and the player navigates each position using local calculation guided by learned intuition. New concepts emerge between games, gradually, through accumulated experience. If we treat each position as a puzzle, the existing vocabulary turns out to be sufficient, at least approximately, and no single position demands a concept that doesn't yet exist.

The same three elements appear in mathematics, but the balance between them is different. A mathematician's conceptual vocabulary and proof techniques are priors, and working through a specific proof attempt (trying constructions, checking special cases, following chains of implication) is local search. The third element, extracting information from search to update understanding, is where the differences between mathematics and chess become most interesting, because in mathematics a single problem can demand a concept that does not yet exist.

The mutilated chessboard problem is a good example of what this looks like.%
\begin{marginfigure}[-10cm]
\centering
\includegraphics[width=\linewidth]{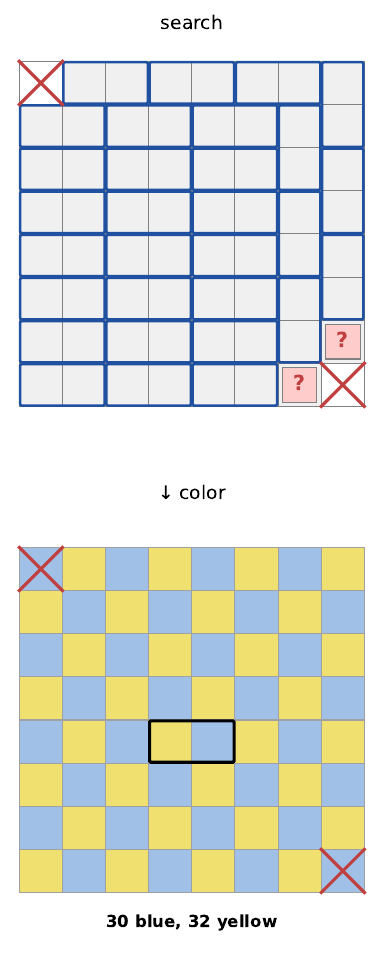}
\caption{The mutilated chessboard. \emph{Top:} a greedy search for a tiling gets stuck with two uncovered squares. \emph{Bottom:} coloring reveals both removed corners are blue. Each domino covers one of each, but $30 \neq 32$.}
\end{marginfigure} Remove two diagonally opposite corners from a chessboard and ask whether the remaining 62 squares can be tiled with dominoes. In the vocabulary of domino placements, each available computation is ``try this particular arrangement,'' and each failed attempt tells you almost nothing about whether any arrangement works. A mathematician might bring a general heuristic to this problem: when a construction seems to keep failing, look for an invariant that explains why. But which invariant? The general heuristic ``look for an invariant'' is a background prior, learned from experience across many problems. The specific invariant, coloring, must be discovered from the structure of this particular problem. Once you color the board, a single computation resolves the question: each domino covers one square of each color, and the two removed corners are the same color, leaving 32 of one color and 30 of the other, so no tiling can exist. The concept of coloring was not available in the vocabulary of placements. It had to be created from the data of this specific problem, guided by the general heuristic but not determined by it. This is the characteristic structure of mathematical problem-solving: general heuristics modulate the search, but the specific concepts needed to make progress must often be learned from the problem at hand.

\section*{A model of mathematical problem-solving}

I want to propose a model of how this process works in practice. What does it mean, concretely, to ``extract information from computation''? There is a puzzle here that is worth thinking through. In Shannon's theory \cite{shannon1948mathematical}, information measures uncertainty in a signal: how surprising it is when you receive it, how many bits you need to encode it. A consequence of this definition is that any deterministic computation performed on data cannot yield more information about that data than was already present, since processing cannot create uncertainty that was not there before. This seems to rule out the idea that computation can be ``informative.'' But when a mathematician says that checking a special case was informative, the uncertainty being resolved is not about the data, which was fully determined, but about a conjecture. The special case she checked was certain; what was uncertain was whether the conjecture holds, and the outcome of the computation provided evidence one way or the other. The information gained is how much this evidence shifts her beliefs, which is something closer to relative or Bayesian information than to Shannon information. A calculation that moves her confidence in a conjecture from 50\% to 95\% has yielded a great deal of information in this sense, even though the computation itself was deterministic and the data it operated on was already fully specified.

I suggest that the mathematician maintains a belief network, a web of confidences over mathematical statements and their logical relationships. Gowers \cite{gowers2023makes} has suggested that these beliefs might be updated in a broadly Bayesian way, though the precise mechanism remains an open question.%
\begin{marginfigure}[-18cm]
\centering
\includegraphics[width=\linewidth]{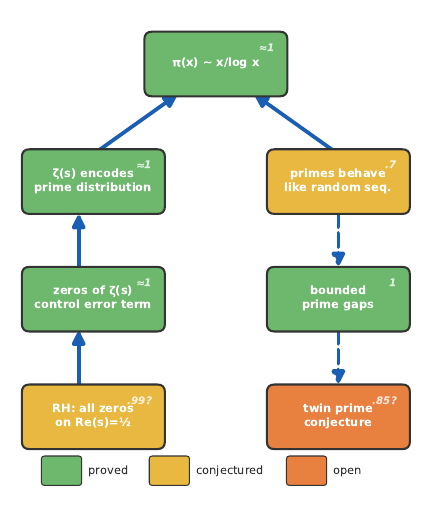}
\caption{A belief network for prime number theory. Arrows show logical implication, but information flows both ways: proving a consequence raises confidence in its premises, and disproving one lowers it.}
\end{marginfigure}
\begin{marginfigure}[-4cm]
\centering
\includegraphics[width=\linewidth]{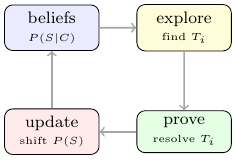}
\caption{The belief-update loop. The mathematician explores for informative auxiliary questions $T_i$, attempts to resolve them, and updates beliefs based on the outcome.}
\end{marginfigure} She acts as an explorer, generating auxiliary questions whose resolution would be informative about the main conjecture: if I check this special case and it works, how much does that increase my confidence? If I try to construct a counterexample and fail, what does the failure tell me? She then acts as a prover, attempting to resolve these auxiliary questions through computation. The outcomes, whether proofs, disproofs, or partial results, update the belief network, and the updated beliefs guide the next round of exploration. This loop continues until the target statement is resolved or until belief drops low enough that she switches to exploring the negation.

How much belief-shifting information each computation yields depends on the vocabulary the mathematician is working in, on which computations are available to her. This is where the framework connects back to the distinction between chess and mathematics: in chess, the loop operates within a fixed vocabulary and each instance is tractable. In mathematics, the loop sometimes requires a vocabulary change to continue making progress.

If this model is roughly correct, it suggests that mathematics is, in an important sense, an experimental science \cite{g2025unreasonableeffectivenessmathematicalexperiments}: the mathematician runs computations (checks special cases, tests conjectures on examples), observes the outcomes, updates beliefs, and uses the updated beliefs to guide the next round of exploration. The ``experiments'' are computations rather than physical measurements, but the logic of the process is the same. And just as in experimental science, the most important advances often come not from running more experiments in the existing framework but from changing the framework itself, from introducing new concepts that make new kinds of experiments possible.

\section*{Concepts and the search landscape}

Why does mathematics require vocabulary changes where chess does not? The fundamental difficulty is exponential blowup: a proof of ten steps where each step offers a hundred possible choices has $100^{10}$ possible paths, far too many for brute force search. Priors help by pruning, but pruning only works within the existing vocabulary. If the computation you need does not exist in your language, no amount of better pruning will find it. You cannot discover that a tiling is impossible by trying placements more cleverly, and no chess-style policy update will teach you the concept of coloring. The solver must create new vocabulary during the process of working on the problem, through local exploration. I will call this explicit reshaping: giving patterns names, creating new vocabulary, and thereby making new computations available.

The K\"onigsberg bridge problem is a case in point. The city of K\"onigsberg had seven bridges connecting four landmasses, and the question of whether one could walk a route crossing each bridge exactly once had attracted many attempts, none successful.

\begin{marginfigure}
\centering
\includegraphics[width=\linewidth]{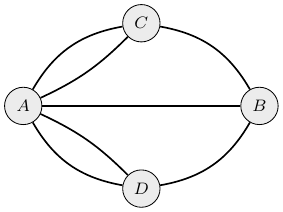}
\caption{The seven bridges of K\"onigsberg. Each vertex has odd degree ($A$: degree 5, $B$, $C$, $D$: degree 3), so no Euler path exists.}
\end{marginfigure}

\noindent In the vocabulary of routes, the available computations are ``try this path and see if it works,'' and each failure provides only the information that this particular route does not work, which says almost nothing about whether any route exists.

Euler's contribution was to notice that the failures had a pattern. Every time you enter a landmass you must also leave it, so each intermediate landmass needs an even number of bridges, and three of K\"onigsberg's four landmasses had an odd number, which meant no route could exist. He extracted this pattern, named it (what we now call the degree of a vertex), and the named pattern made a new computation available: instead of searching for routes, count the vertices with odd degree, and if more than two, no route exists. This single computation, which takes seconds, shifted beliefs about not just K\"onigsberg but every such problem from uncertain to certain. The new vocabulary (graph, vertex, degree) made a class of computations expressible that carried enormously more information per unit of effort than anything the old vocabulary could access.

Current AI systems sit at various points within this framework. AlphaGo reshapes its vocabulary implicitly through training, but within any single game the vocabulary is fixed, and this is sufficient for Go. Large language models have strong priors and can search through chain-of-thought reasoning, but like AlphaGo during a game, their concepts are fixed during any single episode of reasoning. \href{https://deepmind.google/discover/blog/ai-solves-imo-problems-at-silver-medal-level/}{AlphaProof} adds formal verification, giving it a signal for self-play within Lean's proof framework, but the framework itself is fixed. All of these systems operate, in this respect, like a chess player during a game, navigating each instance using existing concepts and local search, and what remains out of reach is the step Euler took: creating a new concept during the process of working on a single problem.

Human mathematicians can take this step. When a mathematician spends a long time on a problem without making progress, working through auxiliary conjectures and special cases and finding that nothing yields a sufficient belief update, the failure itself becomes a signal. It suggests that the available vocabulary does not contain the concepts needed to decompose the problem, and this recognition is what drives the creation of new frameworks.

\section*{The belief-update loop at historical scale}

The belief-update loop often operates at the scale of decades or centuries. As an example, consider the distribution of prime numbers.

Gauss counted primes in large tables and noticed that primes near a large number $X$ appear with frequency roughly $1/\log X$. This was a computation in the vocabulary of counting, and it shifted beliefs: the primes are not random, they thin out in a predictable way. But the vocabulary of counting could not explain \emph{why} this pattern holds, or answer the deeper questions it raises: how close is the approximation? What is the true error term? Where exactly do primes cluster and where do they avoid each other?

\begin{marginfigure}[6cm]
\centering
\includegraphics[width=\linewidth]{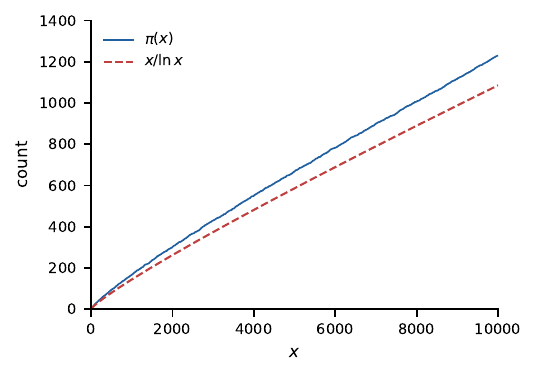}
\caption{The prime counting function $\pi(x)$ (solid) and Gauss's approximation $x/\ln x$ (dashed). The pattern is visible in the data long before the theory explains it.}
\end{marginfigure}

The reformulation of the problem in the language of the Riemann zeta function $\zeta(s) = \sum n^{-s}$ changed what computations were available. In this new vocabulary, questions about the distribution of primes become questions about the zeros of an analytic function: where does $\zeta(s) = 0$? The Riemann hypothesis, which asserts that all nontrivial zeros lie on a specific line, would, if true, give the sharpest possible estimate of the error in Gauss's approximation. The language of the zeta function did not just compress what Gauss had observed. It opened an entirely new set of questions, for instance the value of the zeta function on the critical line. This is a very natural question in the language of complex functions and quite convoluted if expressed directly in terms of the prime numbers - and nevertheless informative of the distribution of prime numbers. Each answer to one of these questions shifts beliefs about the distribution of primes far more efficiently than any computation in the original vocabulary of counting.

An often mentioned criterion for such conceptual understanding is the notion of \textit{Kolmogorov complexity}. There is indeed something to this: a conceptual language and a set of priors that allow for a short search path towards a solution provide a short description of that solution, and thus an upper bound on its Kolmogorov complexity. Conversely, Kolmogorov complexity provides a lower bound on the difficulty of locating a solution through any search process against a fixed prior. However, I am arguing that concepts are good to the extent that they give us good priors, and this is something that Kolmogorov complexity is blind to, along with the related problem of finding the shortest program in the first place. Finally, Kolmogorov complexity is entirely blind to the beauty of the space a set of concepts can break open. A concept that merely compresses is useful, but a concept that opens new questions and reveals shared structure across domains is generative, and it is the generative capacity of a concept, not its compression ratio, that determines how much mathematics it creates.

\section*{What is a concept?}

I have been using the word ``concept'' implicitly so far but to make progress, I need to be more precise. The word ``concept'' has a long and contested history in philosophy. Classical definitions going back to Aristotle treat concepts as defined by necessary and sufficient conditions, which works well for mathematical definitions but poorly for most natural categories. Wittgenstein argued that many concepts, like ``game,'' resist definition entirely and can only be recognized through overlapping family resemblances \cite{wittgenstein1953philosophical}. More recently, prototype theories and theory-based accounts have each captured some aspects of how humans categorize while missing others. In mathematics, Lakatos showed in \emph{Proofs and Refutations} \cite{lakatos1963proofs} that mathematical concepts are not fixed in advance but evolve through the process of proving, shaped by the role they play in arguments. I will not try to resolve these debates, but Lakatos's observation points toward the kind of definition I want: a functional one, motivated by the framework I have been developing. A concept, for my purposes, is anything that reshapes a search landscape.

To make this precise, it helps to think about what a search landscape consists of. Any problem-solving process can be modeled as search in a discrete language of moves. In chess, the language consists of the legal moves available from each position. In mathematics at the lowest level, the language consists of valid proof steps: applying axioms, invoking definitions, using inference rules. At higher levels, the language can include composite moves like ``apply induction on the degree'' or ``try the probabilistic method,'' each of which corresponds to a structured sequence of lower-level steps. The search landscape is the tree of all possible paths through this language, and problem-solving consists of exploring this tree efficiently enough to find a path that reaches the goal.

Within this framework, I will distinguish two forms of concept.

An \emph{implicit concept} changes how you prune within a fixed language. The set of available moves stays the same, but you assign different priorities to different branches, exploring some moves first and skipping others. AlphaGo's policy function is a collection of implicit concepts in this sense: it assigns probabilities to the (unchanged) legal moves, concentrating search on the most promising branches. A chess player's intuition that a position ``looks dangerous'' works similarly, pruning the search tree without adding or removing any legal moves, and AlphaGo's superhuman performance in Go shows that implicit pruning alone can be enormously effective.

An \emph{explicit concept} changes the language itself, introducing new moves that were not previously available. Euler's concept of vertex degree did not improve the search for routes through K\"onigsberg but rather introduced new primitives (graph, vertex, degree) that made a new kind of move available, ``count the vertices with odd degree,'' which did not exist in the route-trying language and required a new language to express. Descartes' coordinates similarly did not improve compass-and-straightedge constructions but replaced them with a different language where ``solve a system of equations'' is an available move, turning open-ended geometric search into routine algebra. An especially accessible example is positional notation for numbers: in Roman numerals, multiplication requires cumbersome manipulation of symbols (try multiplying XLVII by MCXII), while in positional notation, multiplication is a routine computation that a child can learn. More significantly, testing for divisibility by a fixed number, say $3$, is a laborious process in the Roman numeral system and there are no easy patterns to be spotted. In the Arabic numeral system however, divisibility by $3$ can be tested easily by recursively summing up the digits of the number and testing for divisibility. Most importantly, one can imagine spotting this pattern through examples without being aware of it in the first place. The Arabic numeral system makes it easier to search through the landscape of number theoretic patterns. 

These two forms are connected through a natural lifecycle. Euler's introduction of vertex degree changed the language available for reasoning about networks, and at the time this was a genuinely new idea that had to be worked out explicitly. But a modern graph theory student uses vertex degree without thinking about it, the way a chess player uses the concept of a fork. What was once an explicit creation has become an implicit tool, part of the inherited vocabulary, and her intuition tells her when to apply it. This is how mathematical vocabulary grows over time, with each generation absorbing the explicit creations of the last. In fact, we learn implicit concepts mainly through explicit examples of these concepts being used in the intended manner.

On a pedagogical note, one might imagine that we humans learn concepts by being informed of the formal rules and definitions, and we internalize them through a conscious act of memorization. In practice, nothing could be further than the truth. While there is certainly a part for formal, explicit rules in learning, we gain fluency by a subconscious process of familiarity through many examples - in every field from language to abstract mathematics.

What makes a concept good, whether implicit or explicit? Within the information framework I have been developing, a concept is good to the extent that it increases the information gained per unit of search effort, and several properties contribute to this and what follows is certainly a non-exhaustive list.

The first is information density: a single step in the new vocabulary carries as much information as many steps in the old. Counting colors is one computation that resolves the chessboard problem entirely, while trying placements requires one computation per arrangement, each carrying almost no information. Solving a system of equations finds all intersection points in one procedure, while compass-and-straightedge constructions require a new approach for each pair of curves.

The second is breadth: the concept reshapes the landscape for a family of problems rather than just one. Vertex degree works for all Euler-path problems, not just K\"onigsberg. Positional notation works for all arithmetic, not just one multiplication. A concept that helps with only a single problem is closer to a trick, while a concept that helps with a class of problems represents a genuine reshaping of the landscape.

The third is decomposability: the concept breaks a hard problem into checkable sub-problems, each tractable and informative. Before the germ theory of disease, treatment failure was uninformative: ``the patient didn't improve'' gave no indication of why. Once disease was decomposed into specific agents and mechanisms, each step became independently checkable: is the bacterium present? Does this antibiotic kill it in culture? Is the immune response adequate? The naming of intermediate concepts (infection, resistance, immune response) is what made medicine repairable, and the same principle applies in mathematics: a proof organized around named intermediate results has natural checkpoints that make the belief-update loop efficient.

\begin{marginfigure}[0cm]
\centering
\includegraphics[width=\linewidth]{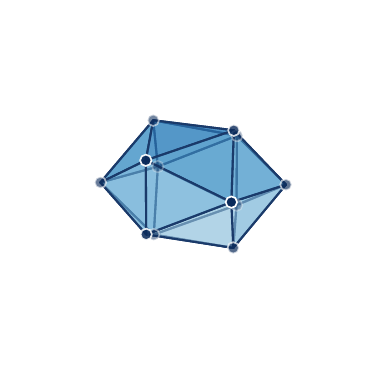}
\caption{An icosahedron drawn on a sphere: $V=12$, $E=30$, $F=20$. The sphere has no holes ($g=0$), and $V-E+F = 2$.}
\end{marginfigure}

The fourth is generativity: the concept opens questions that did not exist in the old language. The concept of a graph, which Euler created to solve the K\"onigsberg bridge problem, is a striking example. Once you have the vocabulary of vertices and edges, questions arise that could not have been asked in the language of bridges and landmasses. How many edges can a graph with $n$ vertices have? When is a graph connected? Can you color the vertices with $k$ colors so that no two adjacent vertices share a color? Each question leads to further structure, and the definitions needed to answer them may themselves be generative. One particularly beautiful question is: if you draw a graph on a surface, dividing it into faces, what is the relationship between the number of vertices $V$, edges $E$, and faces $F$? For any graph drawn on a sphere, the answer turns out to be $V - E + F = 2$, regardless of the graph. A cube has $V=8$, $E=12$, $F=6$, and $8-12+6=2$. A tetrahedron has $V=4$, $E=6$, $F=4$, and $4-6+4=2$. This invariant, discovered by Euler, opened the study of surfaces and eventually algebraic topology, while the same concept of a graph, applied to networks, now powers algorithms like Google's PageRank. The best concepts create far more mathematics than they organize.%

The fifth is transfer: the concept reveals that problems in apparently different fields share structure when expressed in the new vocabulary. Euler's formula $V - E + F = 2$ is a case in point. It began as a fact about graphs drawn on a sphere, a combinatorial observation about counting. But the same formula, applied to a torus (a surface with one hole), gives $V - E + F = 0$. For a surface with $g$ holes, it gives $V - E + F = 2 - 2g$. The number $g$, called the genus, is a property of the surface, not of the graph, and it connects two apparently unrelated fields: combinatorics (counting vertices, edges, and faces) and topology (classifying surfaces by how many holes they have). The formula tells you that if you draw any graph on a surface, the combinatorial count $V-E+F$ reveals the topological shape of the surface. This connection between counting and shape, which could not have been spotted without the vocabulary of graphs, opened the field of algebraic topology and continues to find applications in data analysis and physics.

\begin{marginfigure}[-5cm]
\centering
\includegraphics[width=\linewidth]{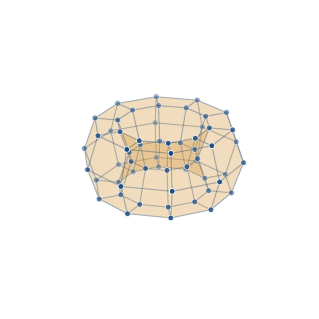}
\vspace{-0.3cm}
\caption{A grid on a torus: $V\!=\!60$, $E\!=\!120$, $F\!=\!60$. One hole ($g\!=\!1$), $V\!-\!E\!+\!F = 0 = 2-2g$.}
\end{marginfigure}

\section*{The value of explicit concepts}

Given that both implicit and explicit concepts can reshape the search landscape, and that implicit concepts alone suffice for domains like chess and Go, it is natural to ask why explicit concepts matter, and I see three connected reasons.

The first is necessity. Some problems require a change in the language itself, not just better pruning within the existing language. No amount of improved intuition about which routes to try would have solved the K\"onigsberg problem, because the information needed to resolve it is not accessible through any computation in the route-trying language. The computation ``count odd-degree vertices'' does not exist in that language and requires new primitives to express. For problems of this kind, implicit reshaping is not enough, no matter how effective, because the required computation is inexpressible in the current vocabulary. I believe that sufficiently complex problems, of the kind found in mathematics and theoretical physics, regularly require this kind of language change, and that this is what distinguishes them from domains like chess and Go where implicit reshaping suffices.

The second is shareability. An explicit concept, because it is named and expressible, can travel between minds. When Euler published his analysis, anyone could use degree-counting to solve bridge-crossing problems, and those who learned the concept could combine it with their own ideas in ways Euler never anticipated. An implicit concept, however effective within one mind or one network's weights, cannot participate in this process. Mathematics advances as a cultural enterprise because concepts are portable across minds, enabling a cumulative accumulation of vocabulary that no single agent could develop alone. This suggests one more factor that contributes to a concept being "good": how well can a concept be internalized from explicit examples.

The third is composability. Explicit concepts can be combined with other explicit concepts to form new ones in a way that implicit concepts cannot, and the resulting combinations can be more powerful than any of their components. Mendeleev's periodic table is a good example from outside mathematics. Chemists already had two named concepts: elements have measurable atomic weights, and certain elements share chemical properties (lithium, sodium, and potassium all react violently with water). Neither concept alone predicted anything new. But Mendeleev composed them, organizing elements into a table where atomic weight determines position and columns share properties. The composed concept had gaps, and those gaps made specific predictions: Mendeleev predicted the existence, atomic weight, density, and chemical behavior of an element he called eka-aluminum. When gallium was discovered in 1875, its properties matched his predictions almost exactly. The composition of two named concepts created predictive power that neither possessed alone, and this is only possible because both concepts were explicit and could be combined in a named framework. The same principle operates throughout mathematics, where definitions routinely compose earlier definitions to produce new structures.%
\begin{marginfigure}[0cm]
\centering
\includegraphics[width=\linewidth]{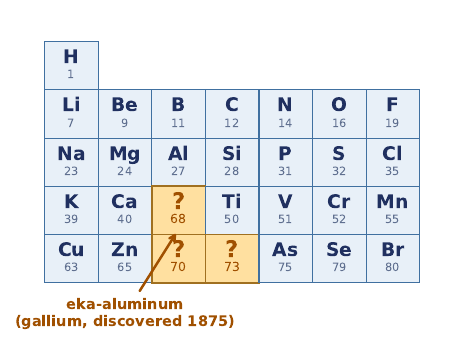}
\caption{Mendeleev's 1869 periodic table. Rows are ordered by atomic weight; columns group elements with similar chemical properties. The gaps predicted undiscovered elements.}
\end{marginfigure}

These three properties are connected in practice. Necessity forces concept creation during the process of solving a problem, which means the concept must be made explicit in order to introduce new computations. Once explicit, it becomes shareable as a byproduct, and once shared, it becomes available for composition with the explicit concepts of others. Individual problems drive explicit concept creation, which feeds cultural accumulation, which expands the vocabulary available for future problems.

\section*{Human and machine mathematics}

Can AI systems learn to create explicit concepts? I can see at least three directions, all research programs rather than working algorithms. The first is automated refactoring: when similar proof structures keep appearing in different contexts, extract the common structure as a definition, much as ``normal subgroup'' likely became a named concept because group theorists kept needing the same hypothesis. The second is introspection: neural networks trained on mathematical proofs develop internal representations that capture structure, and the question is whether these implicit concepts can be made explicit and expressed as definitions, a goal closely related to what mechanistic interpretability research is trying to do. The third is the most ambitious: designing new axiomatic frameworks and evaluating them by how much mathematics they open up, something like what Grothendieck did with schemes.

It is worth noting that the dominant paradigm in AI research today is reinforcement learning, and reinforcement learning is fundamentally the discipline of learning \emph{implicitly} from experience: the agent updates its policy through rewards and penalties, accumulating implicit concepts that reshape its search without ever naming them. The argument of this essay is that for sufficiently complex problems, we also need to think about learning \emph{explicitly} from experience: extracting named patterns from the outcomes of search, creating new vocabulary, and making new computations available. 

Even if machines can learn to do mathematics through the creation of new concepts however, I suspect that humans and machines will approach mathematical problems through very different tradeoffs, and it is worth thinking about what this implies for the future of mathematics.

Humans are severely limited in computational power. Our working memory holds only a handful of items, we process information slowly, and we cannot reliably follow long chains of reasoning. A consequence of this is that we strongly prefer our local searches to be short: we invest in explicit concepts precisely so that the search remaining after reshaping is manageable by a human mind. Human mathematics is concept-heavy precisely because human computation is weak.

The history of computer-assisted proof brings this preference into sharp relief. When Appel and Haken proved the four color theorem in 1976 \cite{appel1977}, the human contribution was a conceptual reduction: they showed, using a technique called discharging, that any counterexample must contain one of a finite set of configurations, and that each such configuration can be shown to be ``reducible'' (meaning it cannot appear in a minimal counterexample). The conceptual insight reduced an infinite problem to a finite one. But the finite problem was still enormous: 1,482 configurations, each requiring detailed verification that took over a thousand hours of computer time. The mathematical community accepted the result only grudgingly, and the dissatisfaction was revealing. The dissatisfaction was not about the computer's accuracy but about the fact that a proof whose key step could not be surveyed by a human mind felt like it was missing something, even if the conclusion was correct. The later simplified proof by Robertson, Sanders, Seymour, and Thomas in 1997 reduced the configuration set to 633, but the essential dissatisfaction remained. Another prominent source of such examples is analytic number theory, where arguments often reduce a problem to checking a finite but very large number of cases. Helfgott's proof \cite{helfgott2014} that every odd number greater than five is the sum of three primes combined the classical circle method with computational verification of all cases below approximately $10^{27}$.

This is of course a spectrum because many proofs involve a reduction to finitely many cases small enough to be checked by the human mind. Whether we find a proof satisfying depends, at least in part, on where it falls on this spectrum. We expend relatively large amounts of energy building enough conceptual scaffolding to reduce the number of cases to a manageable number but this is a subjective threshold. One can therefore imagine that a different intelligent entity, with different computational limits might prefer a different tradeoff, accepting more computational power and less conceptual scaffolding.

Indeed, modern AI systems face the opposite situation. While AI systems are not yet at the point of proving deep mathematical theorems, even the tentative systems we have built often prove theorems by expending much more computational effort than a human would. The situation is even starker in domains like chess and Go, where superhuman play is achieved through a combination of search and pattern recognition that is completely opaque to human understanding. The tradeoff here is that the machine can find correct proofs or winning moves without needing to build the conceptual scaffolding that a human would require, but the resulting arguments are often inscrutable to human minds.

This difference raises a question about what proof is for. Imagine that an AI proves the Riemann Hypothesis with a 10,000-page formal proof in Lean. Every step checks out, and the proof is certified correct. But the argument reveals no structure, introduces no reusable ideas, and offers no insight into why the hypothesis is true. Is this mathematics? If the point of mathematics is to establish which statements are true, the proof is entirely sufficient. If the point is to build understanding, to develop concepts that organize and compress and that make new questions askable, then something essential is missing. The answer depends on what we think mathematics is trying to do.

We stand at another transition point in the history of mathematics, perhaps the most significant since the invention of formal proof. While predictions are always very difficult, they can sometimes serve as targets to move towards, and it is in that spirit that I offer the following two visions of the future.

One is the development of machines that don't merely prove theorems but explain them, extracting from their proofs the patterns and structures that make the results meaningful to us humans and expressing those patterns in a way that we can understand. What counts as a satisfying explanation would be subjective and vary from person to person, and so the machine would have to understand not only mathematics but also the human understanding of mathematics. As a working mathematician, my personal motivation in pursuing mathematics is often to understand things \textit{my way}, and so this vision has at least a modicum of appeal for me. On the other hand, this is far from the only motivation to pursue mathematics, and it is not clear how satisfactory this future would be from many other viewpoints.

The second possibility is that mathematics bifurcates. Machines pursue mathematics, proving ``true'' and ``important'' theorems for applications elsewhere, while humans treat mathematics more as a game, much like chess or Go today when machines are already superhuman. In this future, human mathematics revolves around human contact and the pleasure of understanding mathematics, either individually or in collaboration with others. Perhaps there will be sessions where we put down our computer friends and do mathematics the ``old-fashioned way,'' just as we play chess with friends even though we have superhuman chess engines at our fingertips. This future also has appeal, centered around a different and equally common set of motivations to pursue mathematics.

Of course, it is quite likely that both of these futures come to pass simultaneously, and even more likely that the future is far stranger than what I can imagine. But I think it is worth thinking about these possibilities now, as we build the tools that will shape the future of mathematics and as we decide what to value in mathematical work. The future of mathematics is not just a question of what we can do but also of what we want to do, and the choices we make now will shape the kind of mathematics we have in the future.

\printbibliography

\end{document}